\documentclass{article}%
\usepackage{amsmath}%
\usepackage{amsfonts}%
\usepackage{amssymb}%
\usepackage{graphicx}

\begin{document}

\title{A Note on Laplace Transforms of Some Particular Function Types}
\author{Henrik Stenlund\thanks{The author is obliged to Visilab Signal Technologies for supporting this work.}
\\Visilab Signal Technologies Oy, Finland}
\date{9th February, 2014}
\maketitle
\begin{abstract}
This article handles in a short manner a few Laplace transform pairs and some extensions to the basic equations are developed. They can be applied to a wide variety of functions in order to find the Laplace transform or its inverse when there is a specific type of an implicit function involved.\footnote{Visilab Report \#2014-02}
\subsection{Keywords}
Laplace transform, inverse Laplace transform
\subsection{Mathematical Classification}
MSC: 44A10
\end{abstract}
\section{Introduction}
\subsection{General}
The following two Laplace transforms have appeared in numerous editions and prints of handbooks and textbooks, for decades. 
\begin{equation}
\textsl{L}_{t}[F(t^2)],{s}=\frac{1}{2\sqrt{\pi}}\int^{\infty}_{0}{u^{-\frac{3}{2}}}e^{-\frac{s^2}{4u}}{f(u)}du \ \ (FALSE) \ \label{eqn10}
\end{equation}
\begin{equation}
\textsl{L}^{-1}_{s}[\frac{f(ln(s))}{s\cdot{ln(s)}}],{t}=\int^{\infty}_{0}{\frac{t^{u}{f(u)}du}{\Gamma(u+1)}} \ \ (FALSE) \ \label{eqn20}
\end{equation}
They have mostly been removed from the latest editions and omitted from other new handbooks. However, one fresh edition still carries them \cite{Schaum2013}. Obviously, many people have tried to apply them. No wonder that they are not accepted anymore since they are false. The actual reasons for the errors are not known to the author; possibly it is a misprint inherited from one print to another and then transported to other books, believed to be true. The author tried to apply these transforms, stumbling to a serious conflict. Then the interest arose to sort it out and to be able to use them properly. These transforms are useful in solving integrals and infinite series, in spite of their ugly looks. We mark the transform pairs as usual. We assume all functions to be handled in the following, to fulfill all requirements for the existence of the transform and its inverse, putting aside mathematical rigor.
\begin{equation}
\textsl{L}^{-1}_{s}[f(s)],{t}=F(t) \label{eqn30}
\end{equation}
\begin{equation}
\textsl{L}_{t}[F(t)],{s}=f(s)=\int^{\infty}_{0}{e^{-st}{F(t)}dt}  \label{eqn40}
\end{equation}
A third transform belongs to this category too but seems to be restored to its correct form in the new edition \cite{Schaum2013} of the handbook, the older editions having the false form.
\begin{equation}
\textsl{L}^{-1}_{s}[\frac{f(\sqrt{s})}{s}],{t}=\frac{1}{\sqrt{\pi}{t}}\int^{\infty}_{0}{e^{-\frac{u^2}{4t}}{F(u)}du} \ \ (FALSE) \ \label{eqn50}
\end{equation}
\begin{equation}
\textsl{L}^{-1}_{s}[\frac{f(\sqrt{s})}{\sqrt{s}}],{t}=\frac{1}{\sqrt{\pi}{t}}\int^{\infty}_{0}{e^{-\frac{u^2}{4t}}{F(u)}du} \ \ \ \label{eqn60}
\end{equation}
This transform is not studied further. Its proof is analogous compared to those handled.

In Section 2 we derive the correct forms for equations (\ref{eqn10}) and (\ref{eqn20}). After generalizing these results we exhibit new simple example transforms which may prove useful.
\section{Derivation of Special Laplace Transforms}
\subsection{The Implicit Square Function}
To rectify equation (\ref{eqn10}) we start with a guess for the correct form for the inverse function and progress to the transform function which is not known at this point.
\begin{equation}
\textsl{L}_{t}[{t}F(t^2)],{s}=\int^{\infty}_{0}{e^{-s{t}}{{t}F(t^2)}dt}  \label{eqn100}
\end{equation}
We will change the variable $t$ to $\sqrt{r}$ to get
\begin{equation}
\textsl{L}_{t}[{t}F(t^2)],{s}=\frac{1}{2}\int^{\infty}_{0}{e^{-s\sqrt{r}}{F(r)}dr}  \label{eqn110}
\end{equation}
Then we use the known Laplace transform
\begin{equation}
\textsl{L}_{t}[\frac{a}{2\sqrt{\pi{t^3}}}e^{-\frac{a^2}{4t}}],{s}=e^{-a\sqrt{s}} \ \label{eqn120}
\end{equation}
and place it inside the integral above. We get, after swapping the integrations
\begin{equation}
\textsl{L}_{t}[{t}F(t^2)],{s}=\frac{s}{4\sqrt{\pi}}\int^{\infty}_{0}{e^{-\frac{s^2}{4u}}u^{-\frac{3}{2}}{du}\int^{\infty}_{0}{e^{-uz} {F(z)}dz}}  \label{eqn130}
\end{equation}
The latter integral is obviously the Laplace transform $f(u)$ giving the final result.
\begin{equation}
\textsl{L}_{t}[{t}F(t^2)],{s}=\frac{s}{4\sqrt{\pi}}\int^{\infty}_{0}{e^{-\frac{s^2}{4u}}u^{-\frac{3}{2}}{f(u)}{du}}  \label{eqn140}
\end{equation}
\subsection{The Implicit log Function}
Proving the logarithmic case is not too complicated. We start by subjecting the following expression to the Laplace transform, assuming for a moment that it holds
\begin{equation}
\textsl{L}^{-1}_{s}[\frac{f(ln(s))}{s}],{t}=\int^{\infty}_{0}{\frac{t^{u}{F(u)}du}{\Gamma(u+1)}} \label{eqn200}
\end{equation}
and obtain
\begin{equation}
\frac{f(ln(s))}{s}=\int^{\infty}_{0}{dt\cdot{{e^{-st}}}\int^{\infty}_{0}{\frac{t^{u}{F(u)}du}{\Gamma(u+1)}}} \label{eqn210}
\end{equation}
We swap the integrations
\begin{equation}
\frac{f(ln(s))}{s}=\int^{\infty}_{0}{\frac{{F(u)}du}{\Gamma(u+1)}}\int^{\infty}_{0}{dt\cdot{t^{u}\cdot{{e^{-st}}}}} \label{eqn220}
\end{equation}
and then replace the inner integral's variable to $t=\frac{w}{s}$ getting
\begin{equation}
\frac{f(ln(s))}{s}=\frac{1}{s}\int^{\infty}_{0}{{e^{-u\cdot{ln(s)}}}{F(u)}du} \label{eqn230}
\end{equation}
We have an obvious identity and thus the expression below is the final result.
\begin{equation}
\textsl{L}^{-1}_{s}[\frac{f(ln(s))}{s}],{t}=\int^{\infty}_{0}{\frac{t^{u}{F(u)}du}{\Gamma(u+1)}} \label{eqn240}
\end{equation}
\subsection{Generalizing the Implicit Square Function}
We can add a parameter $a\in C$ to the equation (\ref{eqn140}) and the derivation goes along the same lines to
\begin{equation}
\textsl{L}_{t}[{t}F(a\cdot{t^2})],{s}=\frac{s}{4\sqrt{\pi}\cdot{a^{\frac{3}{2}}}}\int^{\infty}_{0}{e^{-\frac{s^2}{4ua}}u^{-\frac{3}{2}}{f(u)}{du}}  \label{eqn300}
\end{equation}
The parameter can be used for generating more expressions. See the following paragraph for more ways of generalizing. We skip them with this transform since the other transform is more fruitful in this respect.
\subsection{Generalizing the Implicit log Function}
Similarly, we can add a parameter $a\in C$ to equation (\ref{eqn240})
\begin{equation}
\textsl{L}^{-1}_{s}[\frac{f(a\cdot{ln(s)})}{s}],{t}=\int^{\infty}_{0}{\frac{t^{au}{F(u)}du}{\Gamma(au+1)}} \label{eqn400}
\end{equation}
It is useful to add more powers $b\in C, Re(b)\geq{0}$ to the denominator
\begin{equation}
\textsl{L}^{-1}_{s}[\frac{f(a\cdot{ln(s)})}{s^{b+1}}],{t}=\int^{\infty}_{0}{\frac{t^{au+b}{F(u)}du}{\Gamma(au+b+1)}} \label{eqn410}
\end{equation}
We can make a translation $\alpha \in C$ in $s$-space getting
\begin{equation}
\textsl{L}^{-1}_{s}[\frac{f(a\cdot{ln(s-\alpha)})}{(s-\alpha)^{b+1}}],{t}=e^{\alpha\cdot{t}}\int^{\infty}_{0}{\frac{t^{au+b}{F(u)}du}{\Gamma(au+b+1)}} \label{eqn420}
\end{equation}
We can actually create a sequence or an infinite series (if it converges) by adding a number $N+1$ of these expressions, with coefficients $c_n$, to have
\begin{equation}
\textsl{L}^{-1}_{s}[f(a\cdot{ln(s-\alpha)})\sum{_{n=0}^{N}{\frac{c_n}{(s-\alpha)^{n+b+1}}}})],{t}=e^{\alpha\cdot{t}}\sum{_{n=0}^{N}{c_{n}t^n\int^{\infty}_{0}{\frac{t^{au+b}{F(u)}du}{\Gamma(au+b+n+1)}}}} \label{eqn430}
\end{equation}
These changes are either proven in the same way as above or by using the known basic properties of the Laplace transform.
We can differentiate with respect to the parameter $a$ in equation (\ref{eqn410}) getting
\begin{equation}
\textsl{L}^{-1}_{s}[\frac{\cdot{f'(a\cdot{ln(s)})}\cdot{ln(s)}}{s^{b+1}}],{t}=\int^{\infty}_{0}{\frac{du\cdot{{u\cdot{F(u)}\cdot{t^{au+b}}}}}{\Gamma(au+b+1)}{[ln(t)-\frac{\Gamma'(au+b+1)}{\Gamma(au+b+1)}]}} \label{eqn440}
\end{equation}
We can differentiate with respect to the parameter $b$ in equation (\ref{eqn410}) to reach at
\begin{equation}
\textsl{L}^{-1}_{s}[\frac{-ln(s)\cdot{f(a\cdot{ln(s)})}}{s^{b+1}}],{t}=\int^{\infty}_{0}{\frac{du{F(u)\cdot{t^{au+b}}}}{\Gamma(au+b+1)}{[ln(t)-\frac{\Gamma'(au+b+1)}{\Gamma(au+b+1)}]}} \label{eqn450}
\end{equation}
\subsection{Simple Results}
\subsubsection{The Implicit log Function}
To find the inverse Laplace transform of the following, with  $\beta{,\gamma}\in C$
\begin{equation}
\frac{e^{-\beta\cdot{ln(s)}}\cdot{sin(\gamma\cdot{ln(s)})}}{s} \label{eqn500}
\end{equation}
we take as the function $f(s)$
\begin{equation}
f(e^{-\beta\cdot{s}}\cdot{sin(\gamma\cdot{s})}) \label{eqn510}
\end{equation}
We know the inverse of this to be
\begin{equation}
F(t)=\frac{1}{2\cdot{i}}[\delta(t-(\beta-i\gamma))-\delta(t-(\beta+i\gamma))] \label{eqn520}
\end{equation}
We can then apply equation (\ref{eqn240}) to obtain
\begin{equation}
\textsl{L}^{-1}_{s}[\frac{e^{-\beta\cdot{ln(s)}}\cdot{sin(\gamma\cdot{ln(s)})}}{s}],{t}=\frac{1}{2\cdot{i}}[\frac{t^{\beta-i\gamma}}{\Gamma(1+\beta-i\gamma)}-\frac{t^{\beta+i\gamma}}{\Gamma(1+\beta+i\gamma)}] \label{eqn530}
\end{equation}
The natural companion for the above is with 
the corresponding $cos(s)$ function
\begin{equation}
f(e^{-\beta\cdot{s}}\cdot{cos(\gamma\cdot{s})}) \label{eqn535}
\end{equation}
with its inverse transform
\begin{equation}
F(t)=\frac{1}{2}[\delta(t-(\beta-i\gamma))+\delta(t-(\beta+i\gamma))] \label{eqn540}
\end{equation}
and we get the following result
\begin{equation}
\textsl{L}^{-1}_{s}[\frac{e^{-\beta\cdot{ln(s)}}\cdot{cos(\gamma\cdot{ln(s)})}}{s}],{t}=\frac{1}{2}[\frac{t^{\beta-i\gamma}}{\Gamma(1+\beta-i\gamma)}+\frac{t^{\beta+i\gamma}}{\Gamma(1+\beta+i\gamma)}] \label{eqn550}
\end{equation}
We can let the $a$ to approach zero in equation (\ref{eqn440})
\begin{equation}
\textsl{L}^{-1}_{s}[\frac{\cdot{f'(0)}{\cdot{ln(s)}}}{s^{b+1}}],{t}=\int^{\infty}_{0}{\frac{du\cdot{{u\cdot{F(u)}\cdot{t^{b}}}}}{\Gamma(b+1)}{[ln(t)-\frac{\Gamma'(b+1)}{\Gamma(b+1)}]}} \label{eqn560}
\end{equation}
But, since
\begin{equation}
f'(0)=-\int^{\infty}_{0}{du\cdot{u\cdot{F(u)}}} \label{eqn570}
\end{equation}
we get the identity
\begin{equation}
\textsl{L}^{-1}_{s}[\frac{ln(s)}{s^{b+1}}],{t}=-\frac{t^{b}}{\Gamma(b+1)}{[ln(t)-\frac{\Gamma'(b+1)}{\Gamma(b+1)}]} \label{eqn580}
\end{equation}
Since $\Gamma'(1)=-\gamma$ we can obtain a simple result for the special case of $b=0$, already known in tables. 
\begin{equation}
\textsl{L}^{-1}_{s}[\frac{ln(s)}{s}],{t}=-ln(t)-\gamma \label{eqn590}
\end{equation}
The derivatives are known for some other special argument values. Since
\begin{equation}
\Gamma'(s+1)=s\Gamma'(s)+\Gamma(s) \label{eqn595}
\end{equation}
we can derive any other derivative at integer arguments. For example
\begin{equation}
\Gamma'(2)=\Gamma'(1)+\Gamma(1)=-\gamma+1 \label{eqn597}
\end{equation}
giving
\begin{equation}
\textsl{L}^{-1}_{s}[\frac{ln(s)}{s^{2}}],{t}=-\frac{t}{\Gamma(2)}{[ln(t)-\frac{\Gamma'(2)}{\Gamma(2)}]}=-t\cdot{ln(t)}+t-t\gamma \label{eqn599}
\end{equation}
\section{Discussion}
The derivation of the correct forms of the Laplace transforms were found to be straightforward though not obvious. The corrected equations are (\ref{eqn140}) and (\ref{eqn240}). The simple result equation (\ref{eqn580}) is possibly new as also the generalization equations (\ref{eqn300}), (\ref{eqn430}), (\ref{eqn440}) and (\ref{eqn450}).  

\end{document}